\documentclass{article}

\usepackage{theorem}
\usepackage{graphicx}
\usepackage{epsfig}
\usepackage[ansinew]{inputenc}
\usepackage{amssymb}
\usepackage{latexsym}
\usepackage{amsfonts}
\usepackage{amsmath}
\usepackage{color}
\usepackage{colortbl}
\usepackage{graphics}
\usepackage{epsfig}
\usepackage{subfigure}
\usepackage{algorithm}
\usepackage{algorithmic}

\newtheorem{thm}{Theorem}
\newtheorem{prop}{Proposition}
\newtheorem{lem}{Lemma} 
\newtheorem{rmk}{Remark}
\newtheorem{pf}{Proof}

\newtheorem{defin}{Definition}
\newtheorem{coro}{Corollary}

\DeclareMathOperator{\arcsinh}{arcsinh}

\begin{document}
	\thispagestyle{empty}
	\title{Numerical solutions of random mean square Fisher-KPP models with advection}
	\date{}
	\author{{\normalsize
			$\begin{array}{c}
				\mbox{\sc M.-C. Casab\'{a}n,  R. Company, L. J\'{o}dar}\\
				\mbox{Instituto Universitario de Matem\'{a}tica Multidisciplinar} \\
				\mbox{Building 8G Access C 2nd floor} \\
				\mbox{Universitat Polit\`{e}cnica de Val\`{e}ncia} \\
				\mbox{Camino de Vera s/n, 46022 Valencia, Spain} \\
				\mbox{\scriptsize{macabar@imm.upv.es; rcompany@imm.upv.es; ljodar@imm.upv.es}}
			\end{array}$}}
	
	\maketitle
	
	

	

\begin{abstract}
	This paper deals with the construction of numerical stable solutions of random mean square Fisher-KPP models with advection. The construction of the numerical scheme is performed in two stages. Firstly, a semidiscretization technique transforms the original continuous problem into a nonlinear unhomogeneous system of random differential equations. Then, by extending to the random framework the ideas of the exponential time differencing method, a full vector discretization of the problem addresses to a random vector difference scheme.  A sample approach of the random vector difference scheme, the use of properties of Metzler matrices and the logarithmic norm allow the proof of stability of the numerical solutions in the mean square sense. In spite of the computational complexity the results are illustrated by comparing the results with a test problem where the exact solution is known.
	
\end{abstract}


\section{Introduction} \label{sec1}

Time dependent and spatial heterogenous diffusion-advection-reaction models are frequent in genetics and po\-pu\-la\-tion dynamics problems related to emigration or invasion biological problems \cite{Fisher,Malchow2008,Bengfort2016,Okubo}. Branching and species creations although they answer to reorganizations and genetic mutations have a random component.

Most of the results  in population dynamics are treated in the deterministic framework, but the random models deserve attention to picture more realistic situations. In fact, the spatial  heterogeneity, the  contact distribution,  sexual distribution among the population, the rate of increase of the population, the speed of the wind in a direction, and even the carrying capacity, all these issues are essentially random. These facts motivate the study of random partial differential diffusion-advection-reaction models.

Deterministic biological invasions models have been treated both theoretically and numerically in \cite{Fisher,Kolmogorov,Skellam,Aronson1975,Aronson1978,Shigesada,
Weinberger2002,Weinberger2007,Zhao,Bartlett}. Deterministic spatial varying coefficient models linked to heterogeneity have been treated in \cite{Bengfort2016,Parrot,Kinezaki,Jin}. A Numerical Montecarlo approch for these problems has been used in \cite{Faou2009}. Fisher-Kolmogorov-Petrosky-Piskunov (Fisher-KPP) models introducing uncertainty throughout Wiener processes and It\^{o} calculus have been studied in \cite{Bengfort2016,Skellam,Doering,Siekman17,McKean}. Fisher-KPP models with non-Wiener randomness, involving random stationary ergodoc coefficient have been studied teoretically in \cite{Beres_Nadin,Nadin}. In this paper we follow the mean square approach developed for both the ordinary and partial differential case in \cite{JC-Rafa,JC-Braumann,AMM_Xelo_JC_Lucas2016}. This approach has two suitable properties. The first is that our solution coincides with the one of the deterministic case, i.e., when the random data are deterministic. The second, is that if $u_{n}(t)$ is a mean square approximation to the exact mean square solution $u(t)$, then the expectation and the variance of $u_{n}(t)$ converges to the expectation and the variance of $u(t)$, respectively \cite{TTSoong}.

Since early 1960s semidiscretization method, so called the method of lines, has been widely used because it take advantage of the powerful results of systems of ordinary differential equations (ODEs), \cite{chudov}. These advantages come out from both the theoretical and practical points of view, typically studing partial differential equations (PDEs) models, see \cite{Sanz,Calvo,Cox,Hoz} and references therein. The method of lines is easy to use from a practical point of view, i.e. engineering applications, but it presents serious challenges from the analysis point of view because of the fact that when the stepsizes discretization tend to zero, the size of the arrival system of ODEs tends to infinity \cite{DAR}. To our knowledge the method of lines has not be used before in the random numerical analysis framework. One of the aims of this paper is to introduce the method of lines in the framework of PDEs models with randomness, in particular in the Fisher-KPP model. As the best model may be wasted with a disregarded analysis, it is also suitable the numerical analysis in the random model case. 

In this paper we consider a one-dimensional random Fisher-KPP equation modeling the s.p. of the density of population of the invasive species $u(x,t)$ depending on time $t$ and spatial variable $x$
\begin{eqnarray}
u_{t}(x,t) &  = &  D(x)\, u_{xx}(x,t) + B(x)\,u_{x}(x,t) + G(x,u) \,, \quad (x,t)\in{]0,\ell[ \,\times \,]0,+\infty[\,}\,, \label{eq}\\
u(x,0) & = & \Phi_0(x)\,, \quad x\in{[0,\ell]}\,,\label{IC_eq}\\
u(0,t) & = & \Phi_1(t)\,, \quad t\in{]0,+\infty[}\,,\label{CC1eq}\\
u(\ell,t) & = & \Phi_2(t)\,, \quad t\in{]0,+\infty[}\,,
\label{CC2eq}
\end{eqnarray}
where the random reaction term takes the form
\begin{equation} \label{Def_G}
G(x,u)=A(x)\left[\,u(x,t)\,(1-u(x,t))\right]\,.
\end{equation}
The coefficients and initial condition of the random mean square problem \eqref{eq}--\eqref{Def_G} are stochastic processes (s.p's) defined in a complete probability space $(\Omega,\mathcal{F},\mathbb{P})$  and the boundary conditions \eqref{CC1eq}--\eqref{CC2eq} are random variables (r.v.'s). To be more specific, the s.p.'s $D(x)$, $B(x)$, $A(x)$, and $\Phi_0(x)$, are described as continuous functions of $x$ with 1-degree of randomness, that is, every one depends on a single r.v. $\delta(\omega)$, $\rho(\omega)$, $\alpha(\omega)$ and $\gamma(\omega)$, respectively. The same results are available, but with more complicated notation, by considering random functions with a finite degree of randomness, see \cite[p.37]{TTSoong}. In addition, we assume the following conditions over the s.p.'s $D(x)$, $B(x)$, $A(x)$ and $\Phi_{0}(x)$, respectively, satisfy the following conditions on their sample realizations:
 \begin{eqnarray}
 D(x,\omega)=d(x)\, \delta(\omega)\,,  & 0 < d_1 \leq D(x,\omega) \leq d_2 < + \infty\,, & x\in]0,\ell[\,, \\ \label{Conditions_D_CI}
 & &  \text{for almost every (a.e.)}\, \omega\in{\Omega}\,, 
 \nonumber \\
 B(x,\omega)=b(x)\, \rho(\omega)\,,   &  |B(x,\omega)| \leq b_1 < + \infty\,, & x\in]0,\ell[\,, \ \text{for a.e.}\, \omega\in{\Omega}\,,
 \label{Conditions_B_CI}\\
 \nonumber \\
 A(x,\omega)=a(x)\, \alpha(\omega)\,, & 0 \leq a_1 \leq A(x,\omega) \leq a_2 < + \infty \,, & x\in]0,\ell[\,, \ \text{for a.e.}\, \omega\in{\Omega}\,, \label{Conditions_A_CI}\\
\nonumber \\
\Phi_0(x,\omega)=\phi_0(x)\, \gamma(\omega)\,, &  0 \leq  \Phi_0(x,\omega) \leq 1\,, & x\in[0,\ell]\,, \ \text{for a.e.} \omega\in{\Omega}\,. \label{Conditions_Phi0_CI}
 \end{eqnarray}
The random boundary conditions \eqref{CC1eq} and  \eqref{CC2eq}, $\Phi_1(t)$ and $\Phi_2(t)$, are s.p.'s with differentiable realizations $\Phi_i(t,\omega)$ and verifying
\begin{equation} \label{Conditions_CCs}
 0 \leq \Phi_{i}(t,\omega) \leq 1\,, \ \text{for a.e.} \ \omega\in{\Omega}\,, \quad t > 0\,, \quad i=1,2\,.
\end{equation}
Note that for the sake of coherence among the initial and boundary conditions one must satisfy
\begin{equation} \label{Coherence_CCsCI}
\Phi_{0}(0,\omega)=\Phi_{1}(0,\omega)\,, \qquad
\Phi_{0}(\ell,\omega)=\Phi_{2}(0,\omega) \,, \ \text{for a.e.} \ \omega\in{\Omega}\,.
\end{equation}

The organization of the remaining part of the paper is as follows. Section 2 that may be regarded as a preliminary one, includes notations, the statement of previous results and the proof of new results about random systems of ODEs using a mean square approach. In particular, the transformation of a nonlinear inhomogeneous initial value problem for a system of random ODEs into an equivalent random integral equation is shown. Section 3 deals with the semi-discretization, further full discretization, the construction of the numerical scheme and, for the sake of clarity in the presentation, the proof of some technical results that will be used in the next section. The main part of the paper that is Section 4 is devoted to the study of qualitative properties of the numerical solution of the random Fisher-KPP model as well as the random stability of the random numerical solution.

For the sake of clarity in the presentation, we recall some notation, definitions and results that will be used later.

A matrix $M\in{\mathbb{R}^{p \times q}}$ is said to be non-negative, denoted $M \geq \mathrm{O}$, if all its entries $m_{i,j} \geq 0$. A square matrix $M\in{\mathrm{R}^{p \times p}}$ is said to be a Metzler matrix if all non-diagonal entries are non-negative. If $M$ is a Metzler matrix in $\mathbb{R}^{p \times p}$, then its exponential $e^{M\, t} \geq \mathrm{O}$ for all $t \geq 0$, \cite{Kaczorek}. We recall that for a rectangular matrix $M\in{\mathbb{R}^{p \times q}}$, its $\| M \|_{\infty}$ is defined by \begin{equation}  \label{Def_normInf}
\| M \|_{\infty} = \max_{1 \leq i \leq p} \sum_{j=1}^{q} |m_{i,j}|\,,
\end{equation}
see \cite[Chap. 2]{Golup}.
If $M$ is an arbitrary square matrix in $\mathbb{R}^{p \times p}$ its $\mu_{\infty}$-logaritmic norm, denoted by $\mu_{\infty}[M]$, \cite[p.33]{Desoer} can be computed as
\begin{equation} \label{Def_logaritmicNorm}
\mu_{\infty}[M]= \max_{1 \leq i \leq p} \left( m_{i,i} + \sum_{j \neq i} \left| m_{i,j}  \right| \right)\,.
\end{equation}
By \cite{Dahlquist}, we have the inequality
\begin{equation} \label{norm_Inf_Exp}
\left\| e^{M \, t}  \right\|_{\infty} \leq e^{t \, \mu_{\infty}[M]}\,, \quad t \geq 0\,.
\end{equation}
Throughout the paper we will denote $\mathbf{1}_{p}$ the vector in $\mathbb{R}^{p}$ having all its entries equal to one.

\section{Preliminaries and new results on $L_{p}$-random matrix calculus} \label{sec2}
For the sake of clarity in the presentation, we denote $L_p(\Omega)$ the space of all real valued r.v.'s $x: \Omega \rightarrow \mathbb{R}$ of order $p$, endowed with the norm
\begin{equation}  \label{Scalar_norm_p}
 \|x\|_p= \left( \mathbb{E}[|x|^p] \right)^{1/p}=
\left( \int_{\Omega} |x(\omega)|^p  \, f_{x}(\omega)\, \mathrm{d}\omega  \right)^{1/p} < + \infty\,,
\end{equation}
where $\mathbb{E}[\cdot]$ denotes the expectation operator and $f_{x}$ the density function of the r.v. $x$.

The space of all random matrices, $L^{m \times n}_p(\Omega)$ endowed with the norm
\begin{equation}  \label{Def_matrixnorm}
\| X \|_p = \sum_{i=1}^{m}\, \sum_{j=1}^{n} \left\| x_{i,j} \right\|_{p}\,,
\quad
x_{i,j}\in{L_p(\Omega)}\,,
\end{equation}
has a Banach space structure. Although we use the same notation for the  norms $\| \cdot \, \|_p $, no confusion is possible because  lower  case letters are used for scalar quantities  and capital letters are used for matrix quantities.

We recall that the matrix norm $\| \cdot \|_{p}$ in spaces $L_{p}^{m \times n}(\Omega)$ is not submultiplicative
\begin{prop}(\cite[Prop. 1]{AMM_Xelo_JC_Lucas2016}) \label{proposition1_AMM}
Let $X=(x_{i,k})\in{L^{m \times q}_{2p}(\Omega)}$ and $Y=(y_{k,j})\in{L^{q \times n}_{2p}(\Omega)}$. Then
\begin{equation} \label{2.1}
\| X\, Y \|_{p} \leq \| X \|_{2p} \, \| Y \|_{2p}\,.
\end{equation}
\end{prop}

\begin{defin}  \label{Def_continuidad_matrix_sp}
Let $\left\{ X(t)=(x_{i,j}(t)), \, t \in{\mathcal{T}} \right\}$ be a matrix s.p. in $L^{m \times n}_{p}(\Omega)$. We say that  $X(t)$ is $p$-continuous at $t\in{\mathcal{T}}$, if it verifies
\[ \lim_{h \rightarrow 0} \| X(t+h)-X(t) \|_{p}= 0\,, \quad t,\ t+h \in{\mathcal{T}}\,. \]
\end{defin}
Then, it is sufficient that its entries $x_{i,j}(t)$ are all $p$-continuous at $t \in{\mathcal{T}}$ in order to guarantee the $p$-continuity of a matrix s.p. $X(t)$.

\begin{lem} \label{continuidad_producto}
Let $\{P(t), t\in{\mathcal{T}}\}$ and $\{Q(t), t\in{\mathcal{T}}\}$ be matrix s.p.'s in $L^{m \times q}_{2p}(\Omega)$ and $L^{q \times n}_{2p}(\Omega)$, respectively.  Assume that $P(t)$ and $Q(t)$ are $2p$-continuous at $t\in{\mathcal{T}}$. Then, the matrix s.p. $P(t)Q(t)\in{L^{m \times n}_{p}}(\Omega)$ and is a $p$-continuous matrix s.p. at $t\in{\mathcal{T}}$.
\begin{pf}

Using Proposition \ref{proposition1_AMM} one follows
\begin{eqnarray*}
 \| P(t+h)Q(t+h)- P(t)Q(t)  \|_p & = &
 \| P(t+h)(Q(t+h)-Q(t)) + (P(t+h)-P(t))Q(t)  \|_p  \\
 & \leq & \| P(t+h) \|_{2p} \| Q(t+h)-Q(t) \|_{2p}  \\
  & & \quad + \| P(t+h)-P(t) \|_{2p} \| Q(t) \|_{2p}\,,   \ \ t,\ t+h \in{\mathcal{T}}\,,
\end{eqnarray*}

Then, taking into account that $P(t)\in{L^{m \times q}_{2p}}(\Omega)$ and $Q(t)\in{L^{q \times n}_{2p}}(\Omega)$, that is,
 $\| P(t) \|_{2p} < + \infty$ and $\| Q(t) \|_{2p}< + \infty$;
and finally the $2p$-continuity of the matrix s.p.'s $P(t)$ and $Q(t)$, one gets
\[  \lim_{h \rightarrow 0}\| P(t+h)Q(t+h)- P(t)Q(t)  \|_p=0, \quad  t,\ t+h \in{\mathcal{T}}.\quad \square\]
\end{pf}
\end{lem}

The next result is a rule for  differentiability of the product of two differentiable matrix  s.p.'s. and it will play a key role later
\begin{prop}(\cite[Proposition 2]{JC-Braumann} \label{proposition2AMM}
Let $F(t)$ and $G(t)\in{L^{n \times q}_{2p}}(\Omega)$ be  differentiable s.p.'s at $\mathcal{T} \subseteq \mathbb{R}$, being  $F'(t)$ and $G'(t)$ its derivatives, respectively. Then, $H(t)=F(t)G(t)\in{L^{m \times q}_{p}}(\Omega)$ is a  differentiable s.p. and
\[ H'(t)=F'(t)G(t)+F(t)G'(t) \,.\]
\end{prop}

\subsection{Random linear non homogeneous matrix differential systems}  \label{sec2.1}

This subsection deals with the construction of a solution for this type of random linear non homogeneous matrix differential systems
\begin{equation}\label{sistemalineal_noHomogeneo}
\left.
\begin{array}{ccl}
Y'(t) & = & L\, Y(t) +N(t), \\
Y(t_0) & = & Y_0,
\end{array}
\right\}\,,\quad \ t\in{\mathcal{J}(t_0)}=[t_0-\delta,t_0+\delta] \subset \mathbb{R},\ \ t_0>0, \ \delta >0\,,
\end{equation}
where $Y(t)$ and $Y_0$ are in $\in{L^{m \times n}_{p}}(\Omega)$, and the square random matrix $L=(\ell_{i,j})\in{L^{m \times m}_{2p}}(\Omega)$ verifies that for its entries $(\ell_{i,j})$ there exist positive constants $s_{i,j}$, $h_{i,j}$ satisfying that their absolute moments of order $r$ are bounded, that is,
\begin{equation} \label{cotaMomentos_L}
\mathbb{E}\left[ |\ell_{i,j}| ^r \right] \leq s_{i,j} \left( h_{i,j} \right)^r<+\infty \,,\quad \forall
r\geq 0\,, \forall i,j:\, 1\leq i,j \leq m\,.
\end{equation}
Furthermore, assume that matrix s.p. $M(t)\in{L_{2p}^{m \times m}(\Omega)}$ of \eqref{sistemalineal_noHomogeneo} is continua in $\mathcal{J}(t_0)$. Then, under all these hypotheses by \cite[Theorem 1]{AMM_Xelo_JC_Lucas2016}, it is clear that
\begin{equation}  \label{Sol_sistemaHomogeneo}
Y_{1}(t)= e^{L(t-t_0)}Y_0\,,
\end{equation}
is the solution s.p. of the random homogeneous problem
\begin{equation}\label{sistemalineal_Homogeneo}
\left.
\begin{array}{ccl}
Y'(t) & = & L\, Y(t), \\
Y(t_0) & = & Y_0,
\end{array}
\right\}\,,\quad \ t\in{\mathcal{J}(t_0)}\,, \ \ t_0>0\,.
\end{equation}
Now, we show that
\begin{equation}  \label{AuxSol_sistemaNoHomogeneo}
Y_{2}(t)= e^{L(t-t_0)}\int_{t_0}^{t} e^{-L(s-t_0)}\, M(s) \, \textrm{d}s\,,
\end{equation}
is the solution of the non homogeneous problem \eqref{sistemalineal_noHomogeneo}.

 In fact, by \cite[p.~103]{TTSoong}, the m.s. derivative of $Y_{3}(t)=\int_{t_0}^{t} e^{-L(s-t_0)}\, M(s) \, \textrm{d}s$ of \eqref{Sol_sistemaNoHomogeneo}  is
\begin{equation} \label{derivadaY3}
Y_{3}^{\prime}(t)=e^{-L(t-t_0)}\, M(t)\,, \quad t\in{\mathcal{J}(t_0)} \,.
\end{equation}

By \eqref{Sol_sistemaHomogeneo}, \eqref{sistemalineal_Homogeneo}, \eqref{derivadaY3} and Proposition \ref{proposition2AMM}, we have
\begin{eqnarray*}
Y_2^{\prime}(t) & = & e^{L(t-t_0)} \, Y_{3}^{\prime}(t) +
L \, e^{L(t-t_0)} \, Y_{3}(t)\\
& = & e^{L(t-t_0)}\, e^{-L(t-t_0)} \, M(t) + L\, e^{L(t-t_0)}
\int_{t_0}^{t}e^{-L(s-t_0)} \, M(s)\, \mathrm{d}s  \\
& = & M(t) + L\, Y_{2}(t)\,,
\end{eqnarray*}
with $Y_{2}(t_0)=0$.
By superposition and linearity of equation of problem \eqref{sistemalineal_noHomogeneo}, the matrix s.p.
\begin{eqnarray}
\Phi(t) & = & Y_1(t) +Y_2(t) = e^{L(t-t_0)}\, Y_0 +
e^{L(t-t_0)}\, \int_{t_0}^{t} e^{-L(s-t_0)}\, M(s) \, \textrm{d}s \nonumber\\
& = & e^{L(t-t_0)} \, \left\{ Y_0 +  \int_{t_0}^{t} e^{-L(s-t_0)}\, M(s) \, \textrm{d}s \right\}\,,  \label{Sol_sistemaNoHomogeneo}
\end{eqnarray}
is the m.s. solution of the random linear non homogeneous matrix differential system \eqref{sistemalineal_noHomogeneo}, see \cite[Sec.~5.1.2]{TTSoong}.

Summarizing, the following result has been established
\begin{prop}  \label{Prop_Sol_sistemaHomogeneo}
Let  $L\in{L^{m \times m}_{2p}}(\Omega)$, $Y_0\in{L^{m \times n}_{2p}}(\Omega)$ and  $M(t)\in{L^{m \times n}_{2p}}(\Omega)$. Assume that random matrix $L$ verifies condition \eqref{cotaMomentos_L} and the matrix s.p. $M(t)$ is $2p$-continuous in $t\in{\mathcal{J}(t_0)}$. Then, the unique solution of the random linear non homogeneous matrix differential system \eqref{sistemalineal_noHomogeneo} in $L^{m \times n}_{2p}(\Omega)$ is given by \eqref{Sol_sistemaNoHomogeneo}.
\end{prop}

Now, we consider the more general random linear non homogeneous matrix differential system where the right hand side depends on the unknown $Y(t)$
\begin{equation}\label{sistemalinealGeneral_noHomogeneo}
\left.
\begin{array}{ccl}
Y'(t) & = & L\, Y(t) +B\, f(Y(t)), \\
Y(t_0) & = & Y_0,
\end{array}
\right\}\,,\quad \ t\in{\mathcal{J}(t_0)}\,, \ \ t_0>0\,.
\end{equation}
We assume $Y(t)\in{L^{m \times 1}_{p}(\Omega)}$, the random initial condition $Y_0\in{L^{m \times 1}_{2p}(\Omega)}$ and the square random matrices $L,B\in{L^{m \times m}_{2p}(\Omega)}$  are mutually independent with $L$ verifying condition \eqref{cotaMomentos_L}. In \eqref{sistemalinealGeneral_noHomogeneo}, the term  $f(Y(t))$ is a non linear random vector s.p. in $L_{2p}^{m \times 1}(\Omega)$ depending on unknown $Y(t)$.

Inspired by the deterministic method of variation of constants, we look for a solution of \eqref{sistemalinealGeneral_noHomogeneo} of the form
\begin{equation}  \label{CandidataSol_sistemaNoHomogeneo2}
Y(t)=e^{L(t-t_0)}\, C(t)\,,
\end{equation}
being $C(t)\in{L^{m \times 1}_{2p}(\Omega)}$  a differentiable s.p. in $t\in{\mathcal{J}(t_0)}$ with $2p$-derivative $C'(t)$.
Imposing that $Y(t)$, given by \eqref{CandidataSol_sistemaNoHomogeneo2}, is a solution of \eqref{sistemalinealGeneral_noHomogeneo} for the computation of its $p$-derivative $Y'(t)$,
and using Proposition \ref{proposition2AMM}, one gets
\begin{eqnarray*}
Y'(t) &= & L\, e^{L(t-t_0)}\,C(t) + e^{L(t-t_0)}\,C'(t)\,, \\
Y'(t) &= & L\, e^{L(t-t_0)}\,C(t) + B \, f\left(e^{L(t-t_0)}\, C(t)\right)\,,
\end{eqnarray*}
that is $e^{L(t-t_0)}\,C'(t)=B \, f\left(e^{L(t-t_0)}\, C(t)\right)$. Hence

\begin{equation} \label{2.2}
C'(t)= e^{-L(t-t_0)} B \, f\left(e^{L(t-t_0)}\, C(t)\right)\,.
\end{equation}
Assuming that vector s.p. $f\left(e^{L(t-t_0)}\, C(t)\right)$ is $2p$-continua in $\mathcal{J}(t_0)$, we can integrate expression \eqref{2.2} in the interval $[t_0,t]$, $t_0>0$ and then use the initial condition of problem \eqref{sistemalinealGeneral_noHomogeneo}
obtaining
\begin{equation} \label{2.3}
C(t)-Y_0= \int_{t_0}^{t}e^{-L\,(s-t_0)} B \, f\left(Y(s)\right) \mathrm{d}s\,.
\end{equation}
Substituting \eqref{2.3} in \eqref{CandidataSol_sistemaNoHomogeneo2} we have just determined a $p$-solution s.p. of  integral type for the problem \eqref{sistemalinealGeneral_noHomogeneo} given by
\[Y(t)= e^{L(t-t_0)} \left[ Y_0 + \int_{t_0}^{t}
e^{-L\,(s-t_0)} B \, f\left(Y(s) \right)\mathrm{d}s \right] = e^{L(t-t_0)}\, Y_0 + \int_{t_0}^{t} e^{L\,(t-s)} \,B \,
f\left(Y(s) \right)\mathrm{d}s \,.\]

Summarizing, the following result has been established
\begin{prop}  \label{Prop_Sol_sistemaNoHomogeneo}
Let $L\in{L^{m \times m}_{2p}}(\Omega)$, $B\in{L^{m \times m}_{2p}}(\Omega)$ and $Y_0\in{L^{m \times 1}_{2p}}(\Omega)$. Assume that random matrix $L$ verifies condition \eqref{cotaMomentos_L} and the non linear vector s.p. $f(Y(t))\in{L^{m \times 1}_{2p}}(\Omega)$ is $2p$-continuous in $t\in{\mathcal{J}(t_0)}$, $t_0>0$. Then, the random linear non homogeneous matrix differential problem \eqref{sistemalinealGeneral_noHomogeneo} is equivalent
to the non linear random integral equation
\begin{equation} \label{2.4}
Y(t)= e^{L(t-t_0)}\, Y_0 + \int_{t_0}^{t} e^{L\,(t-s)} \,B \, f\left(Y(s) \right)\mathrm{d}s \,, \quad t_0 >0\,.
\end{equation}
\end{prop}

\begin{rmk} \label{Remark1}
Thinking in further applications in section \ref{sec3}, if we consider a partition of the interval $\mathcal{J}(0)=[0,T]$ of the form $t^n=nk$, $n=0, \ldots,N_T$, $k=\frac{T}{N_T}$, we write the solution of the random linear equation of \eqref{sistemalinealGeneral_noHomogeneo} for the time $[t^{n},t^{n+1}]$ and the initial condition at $t^n$
\begin{equation} \label{solRemark1}
Y(t^{n+1})= e^{L(t^{m+1}-t^{n})}  Y(t^n) + \int_{t^n}^{t^{n+1}}
e^{L\,(t^{n+1}-s)} B \, f\left(Y(s)\right) \mathrm{d}s \,.
\end{equation}
Taking into account the substitution $z=t^{n+1}-s$ into \eqref{solRemark1}, one gets
\begin{equation} \label{solFinalRemark1}
Y(t^{n+1})= e^{L\,k}  Y(t^n) + \int_{0}^{k}
e^{-L\,z} B \, f\left(Y(t^{n+1}-z)\right) \mathrm{d}z \,.
\end{equation}
\end{rmk}

\begin{rmk} \label{Remark2}
Taking sample realizations of the random linear integral equation \eqref{solFinalRemark1} for each fixed $\omega\in{\Omega}$, we have
\begin{equation} \label{YRemark2}
Y(t^{n+1},\omega)= e^{L(\omega)\,k}  Y(t^n,\omega) + \int_{0}^{k}
e^{-L(\omega)\,z} B(\omega) \, f\left(Y(t^{n+1}-z,\omega)\right) \mathrm{d}z \,,
\end{equation}
that is, the solution of the sample deterministic problem
\begin{equation} \label{YprimaRemark2}
Y'(t,\omega)=L(\omega)Y(t,\omega) +B(\omega) f(Y(t,\omega))\,, \quad t\in{[t^n,t^{n+1}]}, \ \omega\in{\Omega} \ \text{fixed}.
\end{equation}
for $Y(t^n,\omega)$ given.
\end{rmk}

\section{Discretization and random numerical scheme construction} \label{sec3}
This section deals with the construction of the full discretized random numerical scheme for solving problem \eqref{eq}--\eqref{CC2eq} in two steps. Firstly, we develop a semidiscretization of the PDE into the spatial variable using central difference approximations of the spatial derivatives. Then, the resulting system of random ODE's is discretized in time using a kind of exponential differencing method  \cite{CoxMathews,Pazy}.

Let us consider the uniform partition of the interval $[0,\ell]$, of the form $x_i=ih$, $0 \leq i \leq N$, with $Nh=\ell$. The spatial discretization of a realization $\omega\in{\Omega}$ of equation \eqref{eq} at the mesh points yields
\begin{multline}   \label{Discrete1_Eq}
\frac{du_{i}(t,\omega)}{dt}= \\ D(x_i,\omega)  \frac{u_{i-1}(t,\omega)-2u_{i}(t,\omega)+u_{i+1}(t,\omega)}{h^2}  + B(x_i,\omega)  \frac{u_{i+1}(t,\omega)-u_{i-1}(t,\omega)}{2h}     + G(x_i,u,\omega)\,,
\end{multline}
where $u_i(t,\omega)=u(x_i,t,\omega)$ is the numerical approximation of the solution s.p. of \eqref{eq} for all $t>0$,  and $G(x_i,u,\omega)=A(x_i,\omega)\, u_{i}(t,\omega)(1-u_{i}(t,\omega))$, taking into account \eqref{Def_G}. The resulting semidiscretized system of ODE's in time \eqref{Discrete1_Eq}, can be rewritten in the following vector form, for each fixed $\omega\in{\Omega}$,
\begin{eqnarray}
\left.\begin{array}{lcl}
\dfrac{d}{dt} \mathbf{u}(t,\omega) & = & M(\omega)\, \mathbf{u}(t,\omega) + \mathbf{g}(\mathbf{u},\omega)\,, \qquad t>0\,,\\
\\
\mathbf{u}(0,\omega) & = & \left[ \Phi_0(x_0,\omega), \ldots, \Phi_0(x_N,\omega) \right]^{T}\,,
\end{array} \right\}
\end{eqnarray}
where
\begin{equation} \label{Discrete2_Eq}
\mathbf{u}(t,\omega)= \left[ u_0(t,\omega),\cdots,u_N(t,\omega) \right]^{T}=
\left[ \Phi_1(t,\omega), u_1(t,\omega), \ldots, u_{N-1}(t,\omega), \Phi_2(t,\omega) \right]^{T} \in{\mathbb{R}^{N+1}}\,,
\end{equation}

\begin{eqnarray}  \label{RandomDiscreteMatrix}
M(\omega)=\dfrac{1}{h^2}\,
 \left[ \begin{array}{cccccc}
              0      &  0       &  0      &  0 & \cdots & 0\\
             m_{1,0} &  m_{1,1} & m_{1,2} &  0 & \cdots & 0\\
              0      &  m_{2,1} & m_{2,2} & m_{2,3} & \cdots & 0\\
         \vdots & \vdots & \ddots & \ddots & \ddots & \vdots\\
              0 & 0 & \cdots & m_{N-1,N-2} &  m_{N-1,N-1}  & m_{N-1,N}\\
              0 & 0 & \cdots & 0 & 0 & 0
        \end{array} \right]\in{\mathbb{R}^{(N+1)\times (N+1)}}\,,
\end{eqnarray}
it is a tridiagonal band matrix with two zero rows and random nonzero entries
\begin{eqnarray}  \label{Elements_RandomDiscreteMatrix}
\left.\begin{array}{llcl}
m_{i,i-1}=& m_{i,i-1}(\omega) & = & D(x_i,\omega)-\dfrac{B(x_i,\omega)\, h}{2}\,,\\
m_{i,i} =& m_{i,i}(\omega)  & = & -2D(x_i,\omega) \,,\\
m_{i,i+1} = & m_{i,i+1}(\omega) & = & D(x_i,\omega)+\dfrac{B(x_i,\omega)\, h}{2}\,,\\
\end{array}\right\}\,, \quad 1 \leq i \leq N-1\,,
\end{eqnarray}
and
\begin{equation} \label{RandomDiscreteReactionTerm}
\mathbf{g}(\mathbf{u},\omega)=\mathbf{g}(x,t,\mathbf{u}(t,\omega),\omega)=\left[\Phi_1^{\prime}(t,\omega), G(x_1,u,\omega),\ldots,G(x_{N-1},u,\omega),\Phi_2^{\prime}(t,\omega)\right]^{T}\,,
\end{equation}
where the derivatives $\Phi_i^{\prime}(t,\omega)$ means the classic derivative of the realizations $\Phi_i^{\prime}(\cdot,\omega)$ regarded as functions of $t$.

In order to achieve the full discretized scheme we consider a partition of the time interval $\mathcal{J}(0)=[0,T]$ where $T$ is the target time in the fixed station sence. In accordance with Remark \ref{Remark1}, we denote the time mesh points $t^n=n\,k$, $n=0,\ldots,N_T$, $k=\frac{T}{N_T}$.

Using Remark \ref{Remark2}, with $L(\omega)=M(\omega)$, $\omega\in\Omega$, fixed, $B(\omega)=I$, when $I$ denotes the identity matrix, the sample deterministic differential system
\begin{equation}
\dfrac{d}{dt} \mathbf{u}(t,\omega)  =  M(\omega)\, \mathbf{u}(t,\omega) + \mathbf{g}(\mathbf{u},\omega)\,, \qquad \omega \in{\Omega}\ \text{fixed}, \ t \geq t^{n}\,,
\end{equation}
for a given value $\mathbf{u}(t^n,\omega)$, evaluated at $t=t^{n+1}$ states the deterministic linear integral equation for a fixed $\omega \in{\Omega}$,
\begin{equation} \label{IntegralEq}
\mathbf{u}(t^{n+1},\omega)= e^{M(\omega)\,k} \, \mathbf{u}(t^n,\omega) + \int_{0}^{k}
e^{M(\omega)\,z} \, \mathbf{g}\left(\mathbf{u}(t^{n+1}-z,\omega),\omega \right) \mathrm{d}z \,.
\end{equation}
We approximate $\mathbf{u}(t^{n+1},\omega)$ in \eqref{IntegralEq} substituting the value $\mathbf{g}\left(\mathbf{u}(t^{n+1}-z,\omega),\omega \right)$ for all $z\in[0,k]$ by the corresponding value at $z=k$, $\mathbf{g}\left(\mathbf{u}(t^{n},\omega),\omega\right)$, obtaining the approximation
$\mathbf{u}(t^{n+1},\omega) \approx \mathbf{v}^{n+1}(\omega)$
\begin{equation}  \label{IntegralEq_v}
\mathbf{v}^{n+1}(\omega)= e^{M(\omega)\,k} \, \mathbf{u}(t^n,\omega) + \left( \int_{0}^{k}
e^{M(\omega)\,z} \, \mathrm{d}z \right)\, \mathbf{g}\left(\mathbf{u}(t^{n},\omega),\omega \right)  \,.
\end{equation}
From \cite{CoxMathews}, we have
\begin{equation} \label{IntegralEq_v2}
\mathbf{v}^{n+1}(\omega)=\mathbf{u}(t^{n+1},\omega) + \mathrm{O}(k^2)\,, \qquad \omega\in{\Omega}\ \text{fixed}.
\end{equation}
As $M(\omega)\in{\mathbb{R}^{ (N+1) \times (N+1)}}$ is a singular matrix, instead of computing the exact integral $\int_{0}^{k}
e^{M(\omega)\,z} \, \mathrm{d}z$, we approximate this by using the Simpson quadrature rule
\begin{eqnarray} \label{Cuadratura}
\left. \begin{array}{lcl}
\int_{0}^{k}
e^{M(\omega)\,z} \, \mathrm{d}z & = & k\, \Lambda [M(\omega),k] + \mathrm{O}(k^5)\,,\\
\Lambda [M(\omega),k] & = & \dfrac{1}{6}
\left( I +4\,e^{M(\omega)\,\frac{k}{2}} + e^{M(\omega)\,k} \right)\,,
\end{array} \right\}
\end{eqnarray}
see \cite{Atkinson}. Hence, using \eqref{RandomDiscreteReactionTerm} and \eqref{IntegralEq_v}--\eqref{Cuadratura}, we obtain the following discretization $\mathbf{u}^{n+1}(\omega) \approx \mathbf{u}(t^{n+1},\omega)$
\begin{equation} \label{SampledEsquema}
\mathbf{u}^{n+1}(\omega) = e^{M(\omega)\,k} \ \mathbf{u}^{n}(\omega) \ + \  k\ \Lambda [M(\omega),k]\ \mathbf{g}^{n}(\omega)\,, \qquad 0 \leq n \leq N_{T}-1\,,
\end{equation}
where for a fixed $\omega\in{\Omega}$,
\begin{equation} \label{def_gn}
\mathbf{g}^{n}(\omega)= \left[ \frac{\Delta \Phi_1^n(\omega)}{k}, G(x_1,u_{1}^n,\omega), \cdots, G(x_{N-1},u_{N-1}^n,\omega),\frac{\Delta \Phi_2^n(\omega)}{k} \right]^{T} \in{\mathbb{R}^{N+1}}\,.
\end{equation}
and
\begin{equation}  \label{discretizationDerivative}
\frac{\Delta \Phi_i^n(\omega)}{k} = \frac{\Phi_i(t^{n+1},\omega)-\Phi_i(t^n,\omega)}{k}\,, \quad i=1,2\,,
\end{equation}
denotes the approximation of the derivative $\Phi_i^{\prime}(t,\omega)$.
Note that matrix $M(\omega)$ has its first and last zero rows, the $i-th$ row for $i=0,N$ of matrices $e^{M(\omega)\,k}$ and $\Lambda [M(\omega),k]$ take the form
\begin{equation} \label{e0}
\left( e^{M(\omega)\,k} \right)_0 = \left( \Lambda[M(\omega),k] \right)_0 = [1,0, \ldots, 0]\,,
\end{equation}
and
\begin{equation} \label{eN}
\left( e^{M(\omega)\,k} \right)_N = \left( \Lambda[M(\omega),k] \right)_N = [0,0, \ldots, 1]\,.
\end{equation}
Thurs, first and last rows of $\mathbf{u}^{n+1}(\omega)$ given by \eqref{SampledEsquema} state as
\begin{eqnarray}  \label{u0uNAproximations}
\begin{array}{rcl}
u_0^{n+1}(\omega) & = & u_0^n (\omega) + k\, \dfrac{\Delta \Phi_1^n(\omega)}{k} = u_0^n(\omega) + \Delta \Phi_1^n(\omega) = \Phi_1(t^{n+1},\omega)\,,\\
\\
u_N^{n+1}(\omega) & = & u_N^n (\omega) + k\, \dfrac{\Delta \Phi_2^n(\omega)}{k} = u_N^n(\omega) + \Delta \Phi_2^n(\omega) = \Phi_2(t^{n+1},\omega)\,.
\end{array}
\end{eqnarray}
Recovering the boundary conditions \eqref{CC1eq}--\eqref{CC2eq} at level $t^{n+1}$. This is the reason why we consider the approximation
\eqref{discretizationDerivative} instead of the derivatives $\Phi_i^{\prime}(t^n,\omega)$, $i=1,2$.

The discretization of the initial condition, in agreement with \eqref{IC_eq}--\eqref{CC2eq} and \eqref{Conditions_Phi0_CI}--\eqref{Coherence_CCsCI} takes the form
\begin{equation}
\mathbf{u}^0(\omega)  =  \left[ \Phi_1(0,\omega),\Phi_0(x_1,\omega), \cdots, \Phi_0(N-1,\omega),\Phi_2(0,\omega)
\right]^{T}\,,  \label{CI_discretizada}
\end{equation}
and note that $0 \leq \mathbf{u}^{0}(\omega) \leq \mathbf{1}_{N+1}$ a.e. for $\omega\in{\Omega}$.

For the sake of clarity in the presentation, we include the follwing technical lemma.
\begin{lem}   \label{Lemma1}
Let $M(\omega)$ and $\Lambda [M(\omega),k]$ be the matrices defined by \eqref{RandomDiscreteMatrix}--\eqref{Elements_RandomDiscreteMatrix} and \eqref{Cuadratura} and let us assume that $M(\omega)$ is a Metzler matrix. Then
\begin{itemize}
\item[(i)] $\left( e^{M(\omega)\,k} \right)_{i,r} \geq
\left( e^{M(\omega)\,k/2} \right)_{i,r}\,, \quad r=0,N$.
\item[(ii)] $\left( \Lambda [M(\omega),k]  \right)_{i,r} \leq \left( e^{M(\omega)\,k} \right)_{i,r}\,,
\quad 0 \leq i \leq N; \ r=0,N$.
\end{itemize}
\begin{pf}
From the identity $e^{M(\omega)\,k}=e^{M(\omega)\,k/2}e^{M(\omega)\,k/2}$, it follows that
\begin{eqnarray} \label{proofLemma}
\left( e^{M(\omega)\,k} \right)_{i,r} & = & \sum_{s=0}^{N} \left( e^{M(\omega)\,k/2} \right)_{i,s} \, \left( e^{M(\omega)\,k/2} \right)_{s,r} \nonumber \\  
& = & 
\left( e^{M(\omega)\,k/2} \right)_{i,r} \left( e^{M(\omega)\,k/2} \right)_{r,r} + \sum_{\scriptsize{\begin{array}{l}
s=0 \\
s \neq r
\end{array}}}^{N} \left( e^{M(\omega)\,k/2} \right)_{i,s} +
\left( e^{M(\omega)\,k/2} \right)_{s,r}\,.
\end{eqnarray}
As $M(\omega)$ is a Metzler matrix, from \eqref{e0}, \eqref{eN} and
\eqref{proofLemma} one gets part (i). Note that from \eqref{Cuadratura}
\[ \left( \Lambda [M(\omega),k]  \right)_{0,0} =
\left( \Lambda [M(\omega)\,k]  \right)_{N,N}=
\left( e^{M(\omega)\,k} \right)_{0,0} =1\,, \]
and for $i\neq 0$, $i\neq N$, one gets
\[ \left( \Lambda [M(\omega),k]  \right)_{i,r}=
\frac{1}{6} \left( I + 4\, e^{M(\omega)\,k/2} + e^{M(\omega)\,k}  \right)_{i,r} \leq \frac{5}{6} \left( e^{M(\omega)\,k} \right)_{i,r} <
 \left( e^{M(\omega)\,k} \right)_{i,r} \,.\]
 Hence the result is established.
\end{pf}
\end{lem}

As the use of the exponential of a random matrix requieres the hypothesis \eqref{cotaMomentos_L}, we need to assume this condition in order to transit from the sampled deterministic scheme \eqref{SampledEsquema} to a random scheme. Hence we assume that for each fixed $x\in{[0,\ell]}$, the r.v.'s $D(x)$ and $B(x)$ satisfy \eqref{cotaMomentos_L} and \eqref{SampledEsquema},  obtaining the random difference scheme
\begin{equation} \label{EsquemaRandom}
\mathbf{u}^{n+1}= e^{M\,k} \ \mathbf{u}^{n} \ + \  k\ \Lambda [M,k]\ \mathbf{g}^{n}\,, \qquad 0 \leq n \leq N_T-1\,.
\end{equation}
Here $\mathbf{u}^n$ and $\mathbf{g}^n$ are random vectors in $L_{2p}^{(N+1) \times 1}(\Omega)$ being
\begin{equation} \label{def_gnRandom}
\mathbf{g}^{n}= \left[ \frac{\Delta \Phi_1^n}{k}, G(x_1,u_{1}^n), \cdots, G(x_{N-1},u_{N-1}^n),\frac{\Delta \Phi_2^n}{k} \right]^{T} \,,
\end{equation}
with $G$ is defined by \eqref{Def_G}; $M$ is the random matrix in $L_{2p}^{(N+1) \times (N+1)}(\Omega)$ defined by \eqref{RandomDiscreteMatrix}; and $\Lambda [M,k]$ is the random matrix in $L_{2p}^{(N+1) \times (N+1)}(\Omega)$  given by
\begin{equation}  \label{CuadraturaSimpson}
\Lambda [M,k]  = \dfrac{1}{6}
\left( I +4\,e^{M\,\frac{k}{2}} + e^{M\,k} \right)\,.
\end{equation}
The random initial condition of the random scheme \eqref{EsquemaRandom} becomes
\begin{equation} \label{CI_discretizadaRadom}
\mathbf{u}^0 = \left[ \Phi_1(0,\omega),\Phi_0(x_1,\omega), \cdots, \Phi_0(N-1,\omega),\Phi_2(0,\omega) \right]^{T}\,.
\end{equation}
The numerical solution s.p. of the random problem \eqref{eq}--\eqref{Conditions_CCs} at each mesh point $x_i=ih$ and time $t_{n+1}=(n+1)k$ turns out from \eqref{EsquemaRandom}--\eqref{CI_discretizadaRadom}
\begin{eqnarray} \label{EsquemaRandom_xi}
\left.\begin{array}{rl}
 u_{i}^{n+1}  =  &\left( e^{M\,k}  \right)_{i} \, \mathbf{u}^n + k\, \left( \Lambda [M,k] \right)_{i}\, \mathbf{g}^n\,, \quad  1 \leq i \leq N\,, \  0 \leq n \leq N_T-1\,,\\
 \\
 u_{i}^{0}   =  &\Phi_0(x_i)\,,  \quad 0 \leq i \leq N\,,\\
\end{array} \right\}
\end{eqnarray}
where $\left( e^{M\,k}  \right)_{i}$ and $\left( \Lambda [M,k] \right)_{i}$ denote the $i$-th row of the random matrices $e^{M\,k}$ and $\Lambda [M,k]$, respectively.\\

In the next section, we deal with the numerical analysis of the random scheme \eqref{EsquemaRandom_xi}. For the sake of clarity in the presentation, we recall the concept of $\| \cdot \|_{p}$-stability, in the fixed station sense with respect to the time, of a random numerical scheme. That captures the same essence as in the deterministic case, with the only difference that we consider the $\|  \cdot\|_{p}$ of the random numerical solution.
\begin{defin}  \label{Def_Norm_p_stableNumericalScheme}
The random numerical scheme \eqref{EsquemaRandom_xi} is said to be 
$\| \cdot \|_p$-stable in the fixed station sense in the domain $[0,\ell] \times [0,T]$, if for every partition with $k=\Delta t$, $h=\Delta x$ such that $N_{T}k=T$ and $Nh=\ell$, 
\begin{equation} \label{cotaEstabilidad}
\| u_{i}^{n} \|_{p} \leq C, \qquad 0 \leq i \leq N, \ 0 \leq n \leq N_{T}\,,
\end{equation}
where $C$ is independent of the stepsizes $h$, $k$ and of the time level $n$. 
\end{defin}

\section{Properties of the random numerical scheme} \label{sec4}
In this section we prove that starting from a random initial population $\mathbf{u}^0$, such that $0 \leq \mathbf{u}^0(\omega) \leq \mathbf{1}_{N+1}$ a.e. $\omega\in{\Omega}$, the random sequence solution $\{ \mathbf{u}^n\}_{n \geq 0}$ of the scheme \eqref{EsquemaRandom} satisfies
\begin{equation} \label{condition_un}
0 \leq \mathbf{u}^n(\omega) \leq \mathbf{1}_{N+1}\,, \ \ \text{a.e.} \ \omega\in{\Omega}\,,
\end{equation}
that means, that random population $\mathbf{u}^n$ does not overcome the carrying capacity of the habitat, remaining positive. For the sake of convenience, we will use a sample approach considering the sampled scheme \eqref{SampledEsquema} for a fixed $\omega\in{\Omega}$.

The first objective is to find a sufficient condition on the spatial stepsize, $h$, so that matrix $M(\omega)$ given by \eqref{RandomDiscreteMatrix}--\eqref{Elements_RandomDiscreteMatrix} is a Metzler matrix uniformly a.e. $\omega\in{\Omega}$. From \eqref{RandomDiscreteMatrix} and \eqref{Elements_RandomDiscreteMatrix} note that $M(\omega)$ is a Metzler matrix if
\begin{equation}  \label{ConditionMetzlerMatrix}
\left| B(x_i,\omega)  \right| \frac{h}{2} \leq D(x_i,\omega)\,, \quad 1 \leq i \leq N-1\,.
\end{equation}
Note that for a pure diffusive Fisher-KPP model (i.e. $B(x,\omega)=0$), condition \eqref{ConditionMetzlerMatrix} is satisfied for all $h>0$ due to \eqref{Conditions_D_CI} is assumed. In other case, taking
\begin{equation}  \label{1Condition_h}
0 < h \leq 2 \min \left\{ \frac{D(x,\omega)}{\left| B(x,\omega)  \right|}; \ x\in{]0,\ell[}\,, \ \text{a.e.} \  \omega\in{\Omega} \right\}\,,
\end{equation}
$M(\omega)$ is a Metzler matrix, a.e. $\omega\in{\Omega}$. From hypotheses \eqref{Conditions_B_CI}--\eqref{Conditions_D_CI} this is fullfied if
\begin{equation}  \label{2Condition_h}
0 < h \leq \frac{2\, d_1}{b_1}\,.
\end{equation}
Now, we can observe that $i$-th row of the vector scheme \eqref{SampledEsquema} states as
\begin{equation} \label{SampledEsquema_i-th}
u_{i}^{n+1}(\omega) = \left( e^{M(\omega)\,k}  \right)_{i} \, \mathbf{u}^n(\omega) + k\, \left( \Lambda [M(\omega),k] \right)_{i}\, \mathbf{g}^n(\omega)=\mathcal{G}_{i} \left( u_{0}^{n}(\omega), \cdots, u_{N}^{n}(\omega)\right) \,, \ 0 \leq i \leq N\,.
\end{equation}
The non-linear algebraic system for $u_{j}^n$, $0 \leq j \leq N$, provides the new solution at the time level $n+1$ starting from the $n$-th one. In order to study the rate of change of $\mathcal{G}$ with respect to its arguments, we pay attention to the Jacobian of $\mathcal{G}(\omega)$.
From \eqref{SampledEsquema_i-th} and taking into account \eqref{Def_G} and \eqref{def_gn}, we have for the interior rows, $1 \leq i \leq N-1$, of the Jacobian of $\mathcal{G}(\omega)$
\begin{eqnarray} \label{DerivadasGiOmega_interiorRows}
\frac{\partial{\mathcal{G}_{i}}}{\partial{u_j^n}}(\omega) & = &
\left( e^{M(\omega)\,k}  \right)_{i,j} \, + k \,
  \left( \Lambda [M(\omega),k] \right)_{i,j}\, \frac{\partial{\mathbf{g}_j^n}}{\partial{u_j^n}}(\omega)  \nonumber\\
& = & \left( e^{M(\omega)\,k}  \right)_{i,j} + k \,
  \left( \Lambda [M(\omega),k] \right)_{i,j}\,A(x_j,\omega)(1-2u_{j}^n(\omega))\,, \quad 1 \leq j \leq N-1\,.
\end{eqnarray}
For the columns $j=0$ and $j=N$ of the interior rows, from \eqref{RandomDiscreteMatrix}--\eqref{Elements_RandomDiscreteMatrix} and \eqref{def_gn}, it follows
\begin{equation}  \label{DerivadasGiOmega_contornoColumnas}
\frac{\partial{\mathcal{G}_{i}}}{\partial{u_0^n}}(\omega)=
\left( e^{M(\omega)\,k}  \right)_{i,0}\,, \qquad
\frac{\partial{\mathcal{G}_{i}}}{\partial{u_N^n}}(\omega)=
\left( e^{M(\omega)\,k}  \right)_{i,N}\,, \ \ 1 \leq i \leq N-1\,.
\end{equation}
Finally, taking into account \eqref{u0uNAproximations} we have
\[ u_0^{n+1}(\omega)=u_0^n(\omega)+ \Delta \Phi_1^n(\omega)\,, \quad
u_N^{n+1}(\omega)=u_N^n(\omega)+\Delta \Phi_2^n(\omega)\,, \]
and the rows $i=0$ and $i=N$ of the Jacobian of $\mathcal{G}(\omega)$ take the form
\begin{equation}  \label{DerivadasGiOmega_filasi0iN}
\frac{\partial{\mathcal{G}_{0}}}{\partial{u_j^n}}(\omega)= \delta_{0,j}\,, \qquad
\frac{\partial{\mathcal{G}_{N}}}{\partial{u_j^n}}(\omega)= \delta_{N,j}\,, \ \ 0 \leq j \leq N\,,
\end{equation}
where
\begin{eqnarray*}
\delta_{k,\ell}=
\left\{ \begin{array}{lcl}
1 & \text{if} & k=\ell\,,\\
0 & \text{if} & k\neq \ell\,,
\end{array}\right.
\end{eqnarray*}
is the well-known Kronecker delta.
From \eqref{condition_un}, we have
\[  \left| 1-2 u_{j}^n(\omega) \right|  \leq 1 \,, \]
and thus from \eqref{Conditions_A_CI}, \eqref{Cuadratura}, \eqref{2Condition_h} and \eqref{DerivadasGiOmega_interiorRows} we can write
\begin{equation} \label{cota_DerivadasGi}
\frac{\partial{\mathcal{G}_i}}{\partial{u_j^n}}(\omega) \geq
\left( e^{M(\omega)\,k}  \right)_{i,j} -a_2\, k \,
  \left( \Lambda [M(\omega),k] \right)_{i,j}\,, \qquad 1 \leq i,j \leq N-1\,.
\end{equation}
For the sake of convenience, let us introduce the matrix $\mathcal{H}(k,\omega)\in{\mathbb{R}^{(N+1) \times (N+1)}} $ defined by
\begin{equation}  \label{Def_H}
\mathcal{H}(k,\omega)=  e^{M(\omega)\,k}  -a_2\, k \,  \Lambda [M(\omega),k]\,,
\end{equation}
and note under hypothesis \eqref{2Condition_h} that from \eqref{DerivadasGiOmega_contornoColumnas}, \eqref{DerivadasGiOmega_filasi0iN} and \eqref{cota_DerivadasGi}--\eqref{Def_H} one gets
\begin{equation}  \label{mayoracion_G_H}
\frac{\partial{\mathcal{G}}}{\partial{\mathbf{u}^n}}(\omega) \geq
\mathcal{H}(k,\omega)\,.
\end{equation}
The next objective is to determine the value of the parameter $k$ so that
\begin{equation} \label{Positivity_H}
\mathcal{H}(k,\omega) \geq 0\,, \qquad \text{a.e.} \ \omega\in{\Omega}\,.
\end{equation}
Let $m_0(\omega)$ be defined as
\begin{equation}  \label{m_0}
m_0(\omega)= \frac{1}{h^2} \left\{ \min_{0 \leq i \leq N} (m_{i,i}(\omega)) \right\}\,,
\end{equation}
and, under hypothesis \eqref{2Condition_h}, let $R(\omega) \in{\mathbb{R}^{(N+1) \times (N+1)}}$ be the non-negative matrix defined by
\begin{equation}  \label{def_MatrixR}
R(\omega) =M(\omega) -m_0(\omega)\, I\,.
\end{equation}
Note that from \eqref{def_MatrixR},
\begin{equation}  \label{exponential_R}
e^{M(\omega)\, k} =e^{m_0(\omega)\,k\,I} \, e^{ R(\omega)\,k}
\end{equation}
and hence using \eqref{Cuadratura}, \eqref{Def_H} and \eqref{exponential_R} the Taylor series expansion of $\mathcal{H}(k,\omega)= e^{m_0(\omega)\,k\,I} \, e^{ R(\omega)\,k} -a_2\, k \,  \Lambda [M(\omega),k]$ takes the form
\begin{equation} \label{Taylor}
\mathcal{H}(k,\omega)  =  \mathfrak{h}_0(k)\, I + \sum_{s=1}^{\infty} \mathfrak{h}_s(k) \, \frac{(R(\omega))^s\, k^s}{s!}\,,
\end{equation}
where
\begin{eqnarray}
\mathfrak{h}_0(k) & = & e^{m_0(\omega)\, k} - \frac{a_2\ k}{6} \left( 1+ 4 e^{m_0(\omega)\frac{k}{2}} + e^{m_0(\omega)\, k} \right)\,, \quad s=0\,, \label{TaylorTerm_s0}\\
\mathfrak{h}_s(k) & = & e^{m_0(\omega)\, k} - \frac{a_2\ k}{6} \left( \frac{4}{2^s}+e^{m_0(\omega)\frac{k}{2}} + e^{m_0(\omega)\, k} \right)\,, \quad s \geq 1\,. \label{TaylorTerm_s1sInf}
\end{eqnarray}
Taylor expansion \eqref{TaylorTerm_s0} shows that for some $\epsilon (\omega)$, such that $0 < \epsilon (\omega) < k$,
\begin{equation} \label{2TaylorTerm_s0}
\mathfrak{h}_0(k)  =  1- k\, \left( a_2-m_0(\omega) \right) + k^2 \frac{\mathfrak{h}_0^{''}(\epsilon (\omega))}{2}\,,
\end{equation}
where
\begin{multline}  \label{3TaylorTerm_s0}
\mathfrak{h}_0^{''}(\epsilon (\omega)) = m_0^2(\omega) \, e^{m_0(\omega)\,\epsilon (\omega)} + \frac{a_2\, |m_0(\omega) |}{3} \, e^{m_0(\omega)\,\frac{\epsilon (\omega)}{2}} \\
+ \frac{a_2\, |m_0(\omega)|}{6} \, \left( 2- |m_0(\omega)| \epsilon (\omega) \right)\, \left( e^{m_0(\omega)\,\epsilon (\omega)} +
e^{m_0(\omega)\,\frac{\epsilon (\omega)}{2}} \right)\,.
\end{multline}
Note that the sum of the two first terms of the Taylor expansion of $\mathfrak{h}_0(k) $, given in \eqref{2TaylorTerm_s0},
\begin{equation} \label{1Condition_k}
1 - k\, \left( a_2-m_0(\omega) \right) >0\,, \quad \text{if} \ \
k < \frac{1}{a_2+|m_0(\omega)|}\,,
\end{equation}
and by \eqref{m_0} and \eqref{Conditions_D_CI}, the condition \eqref{1Condition_k} occurs when
\begin{equation}  \label{2Condition_k}
k < \frac{h^2}{2\,d_2+a_2\,h^2}\,.
\end{equation}
Since from \eqref{1Condition_k} the term $(2- |m_0(\omega)| \epsilon (\omega))$ of \eqref{3TaylorTerm_s0} is positive,  from \eqref{2TaylorTerm_s0} and \eqref{3TaylorTerm_s0}, we obtain that $\mathfrak{h}_0(k) \geq 0$. Also from \eqref{TaylorTerm_s0} and \eqref{TaylorTerm_s1sInf} one gets
\begin{equation}  \label{mayoracion_H}
\mathfrak{h}_s(k) \geq  \mathfrak{h}_0(k)  \geq 0\,, \quad s \geq 1\,.
\end{equation}
Summarizing, if the spatial stepsize $h$ satisfies \eqref{2Condition_h} and the time stepsize $k$ satisfies \eqref{2Condition_k}, from \eqref{mayoracion_G_H}, \eqref{Taylor} and \eqref{mayoracion_H} we have
\begin{equation}
\frac{\partial{\mathcal{G}}}{\partial{\mathbf{u}^n}}(\omega) \geq
\mathcal{H}(k,\omega) \geq \mathrm{O}\,, \quad 0 \leq n \leq N_T-1\,.
\end{equation}
Under these conditions, the scalar functions $\mathcal{G}_{i} \left( u_{0}^{n}(\omega), \cdots, u_{N}^{n}(\omega)\right)$ are monotone increasing in each argument $u_j^n(\omega)\in{[0,1]}$, where $0 \leq j \leq N$. Hence, from \eqref{SampledEsquema_i-th} one gets
\begin{equation}
\mathcal{G}_{i} \left(u_0^n,0,\cdots,0,u_N^n \right) \leq u_i^{n+1}(\omega) \leq \mathcal{G}_{i} \left(u_0^n,1,\cdots,1,u_N^n \right)\,, \quad 0\leq i \leq N\,.
\end{equation}
From \eqref{Def_G}, \eqref{Cuadratura}--\eqref{discretizationDerivative} and \eqref{SampledEsquema_i-th} one gets
\begin{eqnarray}  \label{G_i_1}
\begin{array}{rcl}
\mathcal{G}_{i} \left(u_0^n,0,\cdots,0,u_N^n \right) & = & \left(e^{M(\omega)\,k}\right)_{i,0}\, u_0^n +
\left(e^{M(\omega)\,k}\right)_{i,N}\, u_N^n +
k \, \left( \Lambda[M(\omega),k] \right)_{i,0} \dfrac{\Delta \Phi_1^n(\omega)}{k} \\
& & + k \, \left( \Lambda[M(\omega),k] \right)_{i,N} \dfrac{\Delta \Phi_2^n(\omega)}{k} \\
\\
& = & \left( e^{M(\omega)\,k} - \Lambda[M(\omega),k]  \right)_{i,0}\,
u_0^n + \left( \Lambda[M(\omega),k] \right)_{i,0}\,
\Phi_1(t^{n+1},\omega) \\
&  & \qquad
+ \left( e^{M(\omega)\,k} - \Lambda[M(\omega),k]  \right)_{i,N}\,
u_N^n + \left( \Lambda[M(\omega),k] \right)_{i,N}\,
\Phi_2(t^{n+1},\omega)
\end{array}
\end{eqnarray}
From Lemma \ref{Lemma1}-(ii) and \eqref{G_i_1} it follows that
\begin{equation}
\mathcal{G}_{i} \left(u_0^n,0,\cdots,0,u_N^n \right) \geq 0\,.
\end{equation}
On the other hand, in order to obtain an upper bound $u_i^{n+1}$, note that from \eqref{SampledEsquema_i-th} and \eqref{G_i_1}
\begin{eqnarray}  \label{G_i_2}
\begin{array}{rcl}
\mathcal{G}_{i} \left(u_0^n,1,\cdots,1,u_N^n \right) & = &
 \left(e^{M(\omega)\,k}\right)_{i} \left(\begin{array}{c}
 u_0^n\\
 1\\
 \vdots\\
 1\\
 u_N^n
 \end{array} \right) + \left( \Lambda[M(\omega),k] \right)_{i,0} \Delta \Phi_1^n(\omega) \\
 & & \quad  + \left( \Lambda[M(\omega),k] \right)_{i,N} \Delta \Phi_2^n(\omega) \\ 
 \\
 & = & \mathcal{G}_{i} \left(u_0^n,0,\cdots,1,u_N^n \right) + \sum_{j=1}^N \left( e^{M(\omega),k} \right)_{i,j}\,.
\end{array}
\end{eqnarray}
Note that if $\Delta \Phi_i^n (\omega) \leq 0$, $i=1,2$, from \eqref{G_i_1} it follows that
\begin{equation}  \label{G_i_3}
\mathcal{G}_{i} \left(u_0^n,0,\cdots,0,u_N^n \right) \leq
\left( e^{M(\omega),k} \right)_{i,0} \Phi_1^n(t^n,\omega)+
\left( e^{M(\omega),k} \right)_{i,N} \Phi_2^n(t^n,\omega)\,.
\end{equation}
In the case that $\Delta \Phi_1^n (\omega) >0$, from Lemma \ref{Lemma1}-(ii), one gets
\begin{eqnarray}  \label{G_i_4}
\left( e^{M(\omega),k} \right)_{i,0}\, u_0^n +
\left( \Lambda[M(\omega),k] \right)_{i,0}\,\Delta \Phi_1^n(\omega) & \leq &
\left( e^{M(\omega),k} \right)_{i,0}\, u_0^n +
\left( e^{M(\omega),k} \right)_{i,0}\,\Delta \Phi_1^n(\omega) \\ \nonumber
& = & \left( e^{M(\omega),k} \right)_{i,0}\,\Phi_1^n(t^{n+1},\omega)\,.
\end{eqnarray}
Analogously if $\Delta \Phi_2^n (\omega) >0$, one gets
\begin{equation}  \label{G_i_5}
\left( e^{M(\omega),k} \right)_{i,N}\, u_N^n +
\left( \Lambda[M(\omega),k] \right)_{i,N}\,\Delta \Phi_2^n(\omega) \leq
\left( e^{M(\omega),k} \right)_{i,N}\,\Phi_2^n(t^{n+1},\omega)\,.
\end{equation}
From \eqref{G_i_2}--\eqref{G_i_5}, whatever the sign of
$\Delta \Phi_i^n (\omega)$, $i=1,2$, be, taking into account that $0 \leq \Phi_i (t,\omega) \leq 1$, $i=1,2$, $0 \leq t \leq T$, one gets
\[ \mathcal{G}_{i} \left(u_0^n,1,\cdots,1,u_N^n \right)  \leq \sum_{j=0}^N \left( e^{M(\omega)\, k} \right)_{i,j}=
\sum_{j=0}^N \left| \left( e^{M(\omega)\, k} \right)_{i,j} \right|\,, \]
and from \eqref{Def_normInf} and \eqref{norm_Inf_Exp}
\begin{equation}  \label{sum_e}
\sum_{j=0}^N \left| \left( e^{M(\omega)\, k} \right)_{i,j} \right| \leq \left\| e^{M(\omega)\,k} \right\|_{\infty} \leq e^{k\, \mu_{\infty}[M(\omega)]}\,.
\end{equation}

Note that from \eqref{Def_logaritmicNorm}, \eqref{RandomDiscreteMatrix}--\eqref{Elements_RandomDiscreteMatrix} under condition \eqref{2Condition_h}, we have
\begin{equation} \label{LogaritmicNorm_M}
\mu_{\infty}[M(\omega)] =0\,,
\end{equation}
and from \eqref{sum_e} one concludes that
\[ \mathcal{G}_{i} \left(u_0^n,1,\cdots,1,u_N^n\right) \leq 1 \,,\]
that is,
\[0 \leq u_i^{n+1} (\omega) \leq 1 \,, \quad 0 \leq i \leq N\,, \ \ \text{a.e.}\ \omega\in{\Omega}\,.\]
Then, \eqref{condition_un} has been proved.\\

Note that as the full discretization of problem \eqref{eq}--\eqref{Conditions_CCs} requieres the use of Proposition \ref{Prop_Sol_sistemaHomogeneo} and Remarks \ref{Remark1} and \ref{Remark2}, coefficients $D(x)$, $B(x)$ and $A(x)$ must be $2p$-continuous s.p.'s. Hence, summarizing, the following result has been established
\begin{thm}  \label{Theorem_Positivity}
Let $D(x)$, $B(x)$ and $A(x)$ be $2p$-continuous s.p.'s satisfying conditions \eqref{eq}--\eqref{Conditions_A_CI}. In addition let us assume that for each fixed $x\in{[0,\ell]}$ the r.v.'s $D(x)$ and $B(x)$ satisfy \eqref{cotaMomentos_L}. Let us assume initial condition s.p. $\Phi_0(x)$ verifies condition \eqref{Conditions_Phi0_CI} and lies in $L_{2^{N_T}p}^{(N+1) \times 1}(\Omega)$, $N_T=T/k$. Let us assume boundary condition s.p.'s $\Phi_i(t)$, $i=1,2$, verify condition \eqref{Conditions_CCs}, lies in $L_{2^{N_T}p}^{(N+1) \times 1}(\Omega)$ and have differentiable realizations $\Phi^{\prime}_i(t,\omega)$, $i=1,2$. Furthermore let us assume that both initial and boundary condition s.p's verify condition  \eqref{Coherence_CCsCI}.

Then, under discretization stepsize conditions \eqref{2Condition_h} and \eqref{2Condition_k}, the numerical solution of the random problem \eqref{eq}--\eqref{Conditions_CCs} constructed by the sampled vector scheme \eqref{SampledEsquema}--\eqref{CI_discretizada} satisfies for each time level $n$, $0 \leq n \leq N_T$,
\begin{equation}  \label{Positivity_uin}
 0 \leq u_i^n (\omega) \leq 1 \,, \quad 0 \leq i \leq N\,, \ \ \text{a.e.}\ \omega\in{\Omega}\,, 
\end{equation}
\end{thm}

\begin{coro}
With the hypotheses of Theorem \ref{Theorem_Positivity}, the random numerical scheme \eqref{EsquemaRandom_xi} is $\| \cdot \|_p$-stable in the fixed station sense.
\end{coro}
\begin{pf}
From \eqref{Positivity_uin} and \eqref{Scalar_norm_p} it follows that 
\[
\| u_{i}^{n} \|_{p} = \left( \int_{\Omega} \left|u_i^n(\omega) \right|^{p} f_{u_{i}^{n}}(\omega) \, \mathrm{d} \omega \right)^{1/p} \leq 
\left(\int_{\Omega} f_{u_{i}^{n}}(\omega) \, \mathrm{d} \omega \right)^{1/p}=1\,,
\]
and consequently, random numerical scheme \eqref{EsquemaRandom_xi} is $\| \cdot \|_p$-stable with constant $C$ of \eqref{cotaEstabilidad} having the value $C=1$.
\end{pf}

Algorithm \ref{algorithm} summarizes the procedure to compute the expectation and the standard deviation of the stable solution s.p. of random scheme \eqref{EsquemaRandom_xi}.
\\

\begin{algorithm}
\caption{Calculation procedure for expectation and standard deviation of solution s.p. $u_i^n$ of random scheme \eqref{EsquemaRandom_xi}}
\label{algorithm}
\begin{algorithmic}[1]
\REQUIRE Hypotheses of Theorem \ref{Theorem_Positivity}
\STATE Obtain the bounds $a_2$, $b_1$, $d_1$ and $d_2$ in conditions \eqref{Conditions_D_CI}--\eqref{Conditions_A_CI} for $D(x)$, $B(x)$ and $A(x)$;
\STATE Select a spatial stepsize $h=\Delta x$ verifying condition \eqref{2Condition_h};
\STATE Consider a partition of the spatial domain $[0,\ell]$ of the form $x_i=i\,h$, $i=0, \ldots,N$, where the integer $N=\frac{\ell}{h}$ is the number of discrete points in $[0,\ell]$;
\FOR{$i=0$ to $N$}
\STATE  Evaluations of $\Phi_0(x_i)$;
\ENDFOR
\STATE Select a temporal stepsize $k=\Delta t$ verifying condition  \eqref{2Condition_k};
\STATE Choose a time $T$;
\STATE Consider a partition of the temporal interval $[0,T]$ of the form $t^n=n\,k$, $n=0, \ldots,N_{T}$, where the integer $N_T=\frac{T}{k}$ is the number of levels necessary to achieve the time $T$;
\STATE Construct the random tridiagonal band matrix $M$, defined in \eqref{RandomDiscreteMatrix}--\eqref{Elements_RandomDiscreteMatrix};
\STATE Compute $e^{Mk}$;
\STATE Construct the random matrix $\Lambda[M,k]$, given by \eqref{CuadraturaSimpson};
\FOR{$n=1$}
\STATE Construct the initial vector solution s.p.'s $\mathbf{u}^0$, given by \eqref{CI_discretizadaRadom}, from steps 4-6;
\STATE Construct the random reaction term $G(x_i,u_i^0)$, $1 \leq i \leq N$, defined in \eqref{Def_G}, using step 14;
\STATE Compute the boundary condition s.p.'s $\Phi_1(t^1)$ and $\Phi_2(t^1)$;
\STATE Compute the increments $\dfrac{\Delta \Phi_i^0}{k}$, $i=1,2$, defined in \eqref{def_gn} using step 16;
\STATE Construct the random vector $\mathbf{g}^0$, defined in \eqref{def_gnRandom}, using steps 15-17;
\STATE Compute $\mathbf{u}^1$ using the random scheme \eqref{EsquemaRandom}, that is, steps 11, 12, 14 and 18;
\ENDFOR
\FOR{$n=2$ to $n=N_T$}
\STATE Compute the boundary condition s.p.'s $\Phi_1(t^{n-1})$ and $\Phi_2(t^{n-1})$;
\STATE Compute  $\dfrac{\Delta \Phi_i^{n-1}}{k}$, $i=1,2$,  using step 22;
\STATE Take the previously calculated value of $\mathbf{u}^{n-1}$;
\STATE Construct the random reaction term $G(x_i,u_i^{n-1})$, $1 \leq i \leq N$; 
\STATE Construct the random vector $\mathbf{g}^{n-1}$, using steps 23-25;
\STATE Compute $\mathbf{u}^{n}$ using the random scheme \eqref{EsquemaRandom}, that is, steps 11, 12, 24 and 26;
\ENDFOR
\STATE Compute $\mathbb{E}\left[u_i^{N_T}\right]$, $1 \leq i \leq N$;
\STATE Compute $\mathbb{E}\left[\left(u_i^{N_T}\right)^2\right]$, $1 \leq i \leq N$;
\STATE Compute $\sqrt{\mathrm{Var}\left[u_i^{N_T}\right]}=+\sqrt{\mathbb{E}\left[\left(u_i^{N_T}\right)^2\right]-\left( \mathbb{E}\left[u_i^{N_T}\right] \right)^2}$, $1 \leq i \leq N$;

\end{algorithmic}
\end{algorithm}
%
%
\section{Numerical example} \label{sec5}
In order to illustrate the analytic results previously established, we consider the following random test heterogeneous problem with available exact solution.
\begin{eqnarray}
u_{t}(x,t) &  = &  (1+x^2)\, u_{xx}(x,t) + x\,u_{x}(x,t) + G(x,u) \,, \quad (x,t)\in{]0,1[ \,\times \,]0,+\infty[\,}\,, \label{eqEj}\\
u(x,0) & = & \left( 1+e^{\sqrt{a/6}\,\arcsinh(x)}  \right)^{-2}\,, \quad x\in{[0,1]}\,,\label{IC_eqEj}\\
u(0,t) & = & \left( 1+e^{-5at/6} \right)^{-2}\,, \qquad \quad \ \ t\in{]0,+\infty[}\,,\label{CC1eqEj}\\
u(1,t) & = & \left( 1+e^{-5at/6+\sqrt{a/6}\,\arcsinh(1)} \right)^{-2}\,, \quad t\in{]0,+\infty[}\,.
\label{CC2eqEj}
\end{eqnarray}
The problem \eqref{eqEj}--\eqref{CC2eqEj} corresponds with the general problem \eqref{eq}--\eqref{Def_G} considering deterministic functions $D(x)=1+x^2$ and $B(x)=x$, and the r.v. $A(x)=a$. The exact solution s.p. of problem \eqref{eqEj}--\eqref{CC2eqEj} when $a$ is deterministic, see subsection 1.6.8.1 of \cite{Polyanin}, is given by
\begin{equation}  \label{exactSolEj}
u(x,t)= \left( 1+e^{-5at/6+\sqrt{a/6}\,\arcsinh(x)} \right)^{-2}\,, \qquad (x,t)\in{[0,1]\times[0,+\infty[}\,.
\end{equation}
We consider that r.v. $a$ is a Gaussian distribution of mean $\mu=0.75$ and standard deviation $\sigma=0.08$ truncated on the interval $[0.01, 1]$, i.e, $a \sim \mathrm{N}_{[0.01,1]}(0.75,0.08)$. Then it is guaranteed that coefficient $A(x)=a$ is $2p$-continuous.
The initial condition s.p. \eqref{IC_eqEj} and the boundary condition s.p.'s \eqref{CC1eqEj}--\eqref{CC2eqEj} verify \eqref{Conditions_Phi0_CI} and \eqref{Conditions_CCs}, respectively.\\
For the calculations, stepsizes $h=0.1$ and $k=0.002$ have been used in the partition of the domain $[0,1] \,\times \,[0,+\infty[$. Then conditions \eqref{2Condition_h} and \eqref{2Condition_k} are guaranteed, being in our case $a_2=1$ and $b_1=d_1=d_2=1$. The fixed time considered was $T=0.01$. Hence, the number of spatial and temporal subintervals were $N=1/h=10$ and $N_T=T/k=5$, respectively.
The initial and boundary condition s.p.'s \eqref{IC_eqEj}--\eqref{CC2eqEj} lie in $L_{2^{5}p}^{11 \times 1}(\Omega)$ because $a$ is a truncated r.v. In addition, it is guaranteed that boundary conditions \eqref{CC1eqEj}--\eqref{CC2eqEj} have differentiable realizations, and deterministic coefficients $D(x)$ and $B(x)$ are $2p$-continuous because they are continuous. Then the hypotheses of Theorem \eqref{Theorem_Positivity} are satisfied.

\begin{figure}[!ht]
\centering \subfigure[]%
    {\includegraphics[width=.5\textwidth]{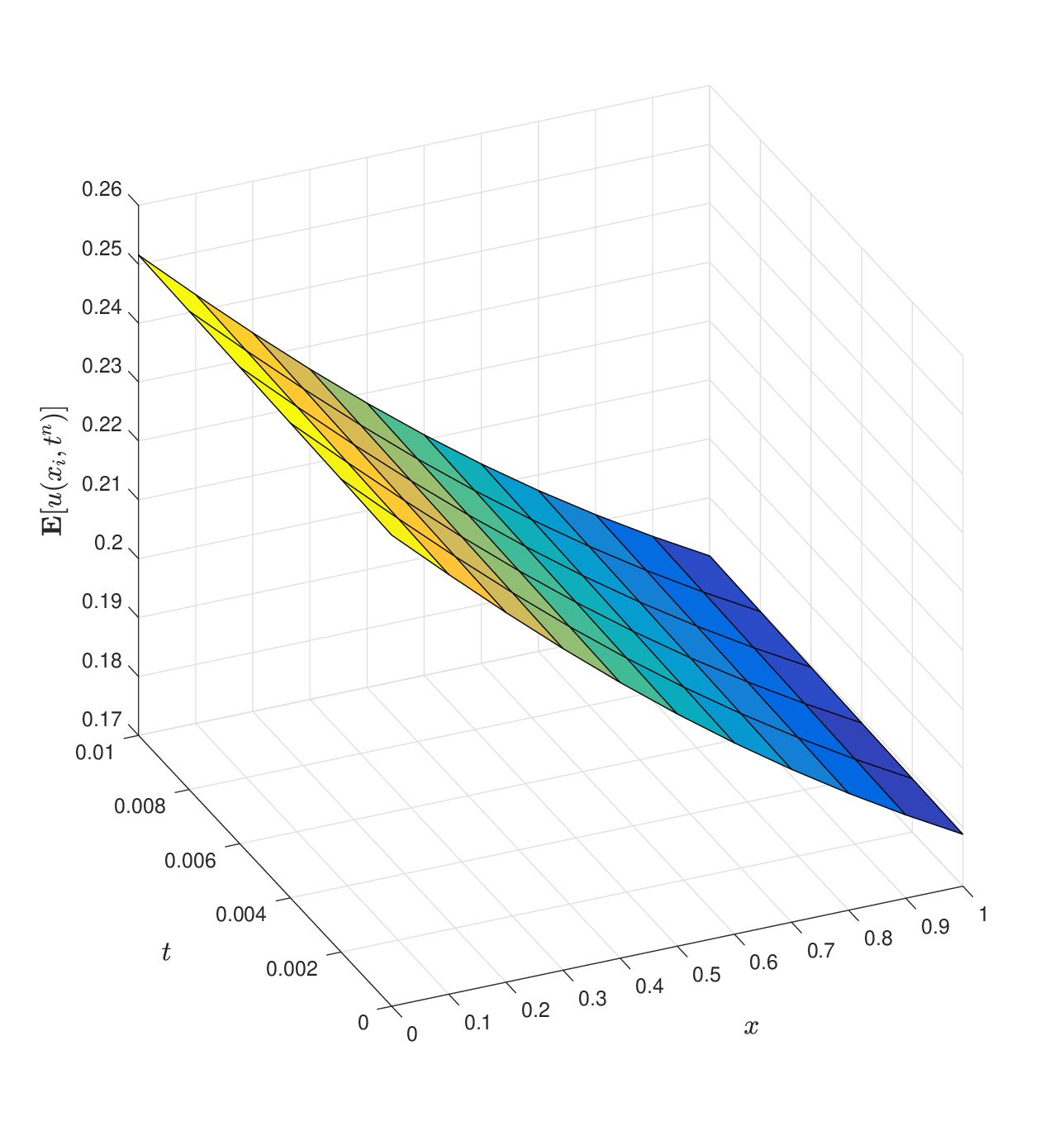}}
\!\!\!\! 
\subfigure[]
{\includegraphics[width=.5\textwidth]{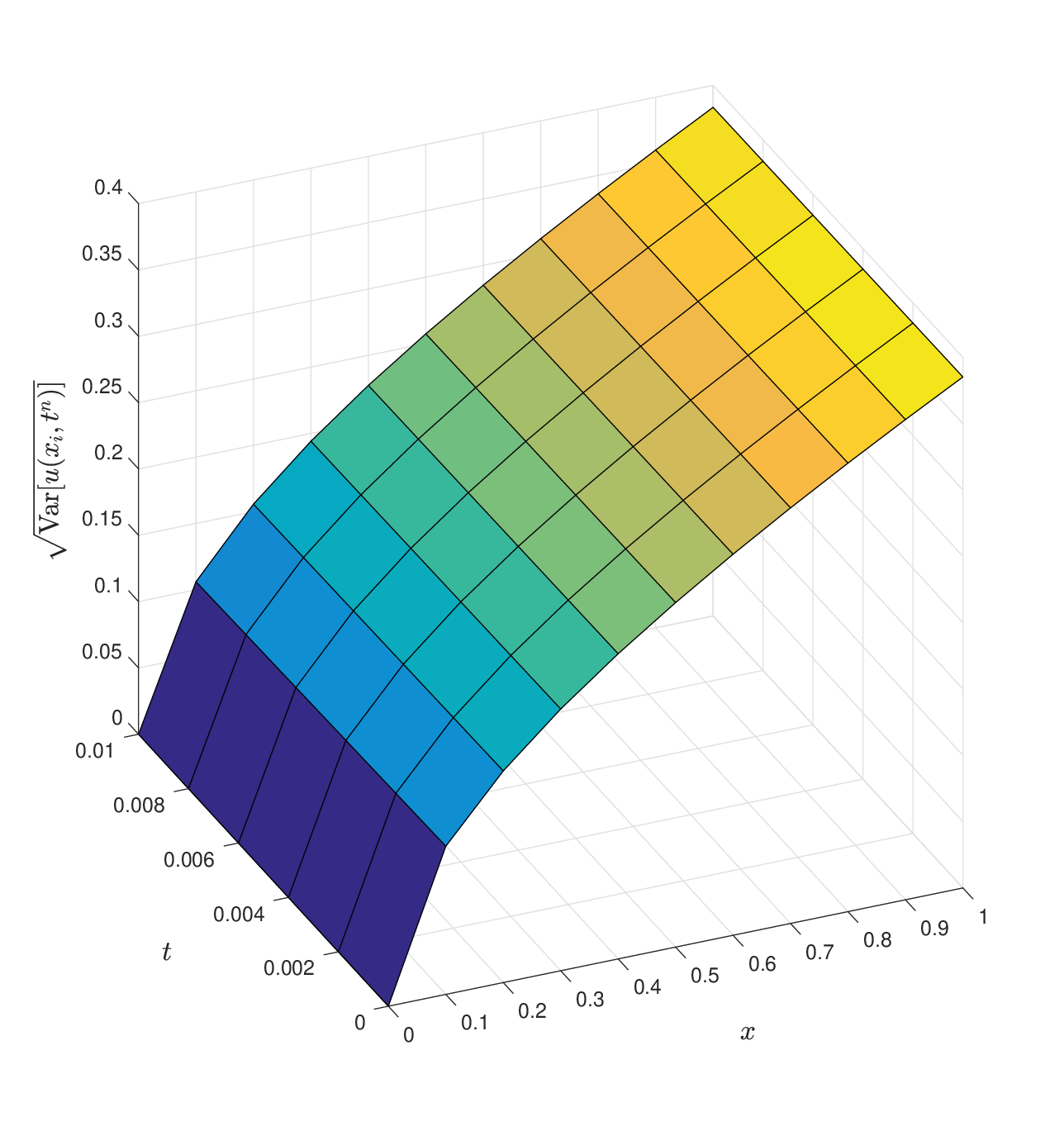}}
 \caption{\label{fig1Example} Plot (a): Surface of the expectation, $\mathbb{E}[u(x_i,t^n)]$. Plot (b): Surface of the standard deviation, $\sqrt{\mathrm{Var}[u(x_i,t^n)]}$. Both
statistical moment functions correspond to the exact solution s.p. \eqref{exactSolEj} of the problem \eqref{eqEj}--\eqref{CC2eqEj}, on the domain $(x_i=ih,t^n=nk)\in{[0,\,1]\times [0, \, T=0.01]}$ for $0 \leq i \leq N=10$ and $0 \leq n \leq N_T=5$.}
\end{figure}

\begin{figure}[!ht]
\centering \subfigure[]%
    {\includegraphics[width=.5\textwidth]{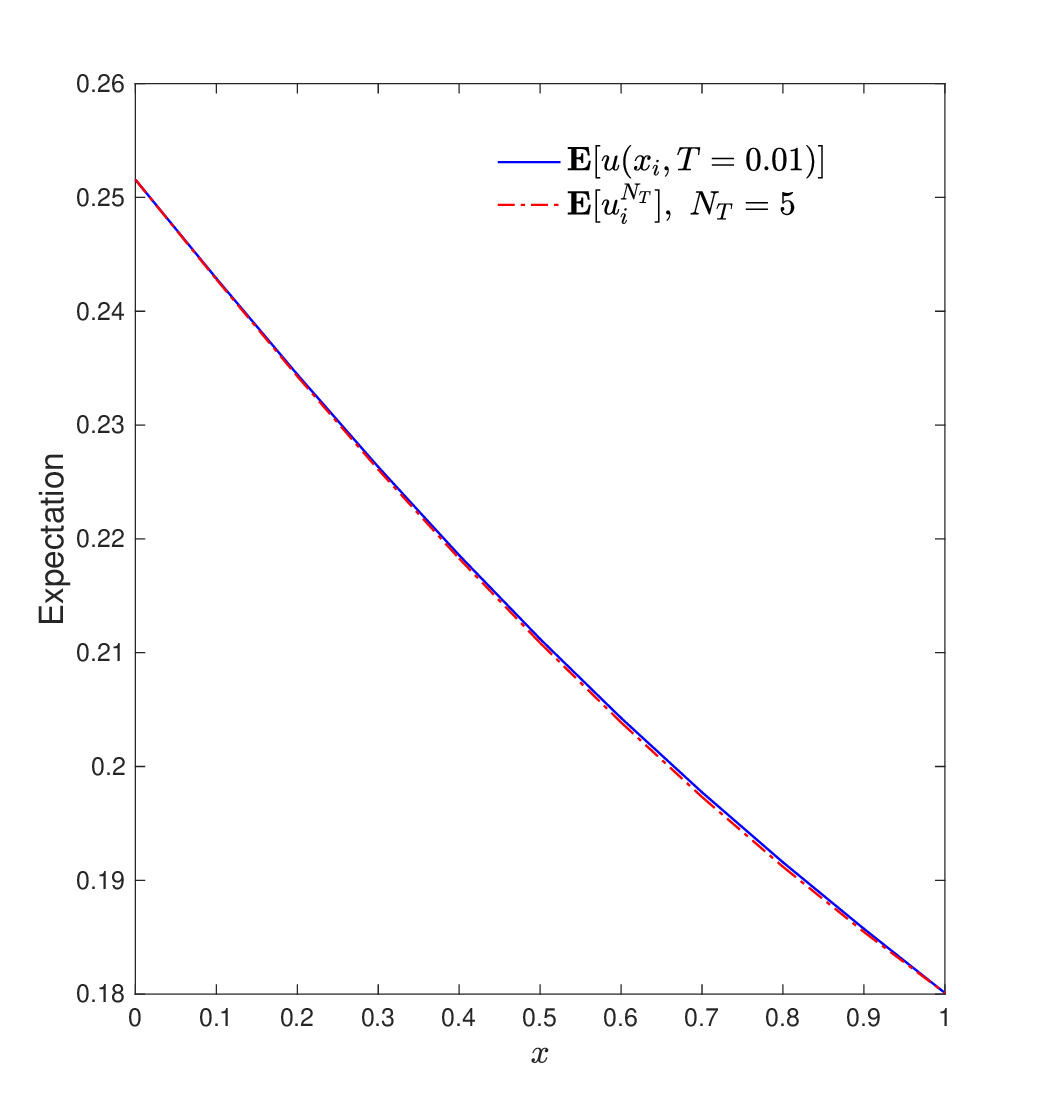}}
\!\!\!\! 
\subfigure[]
{\includegraphics[width=.5\textwidth]
{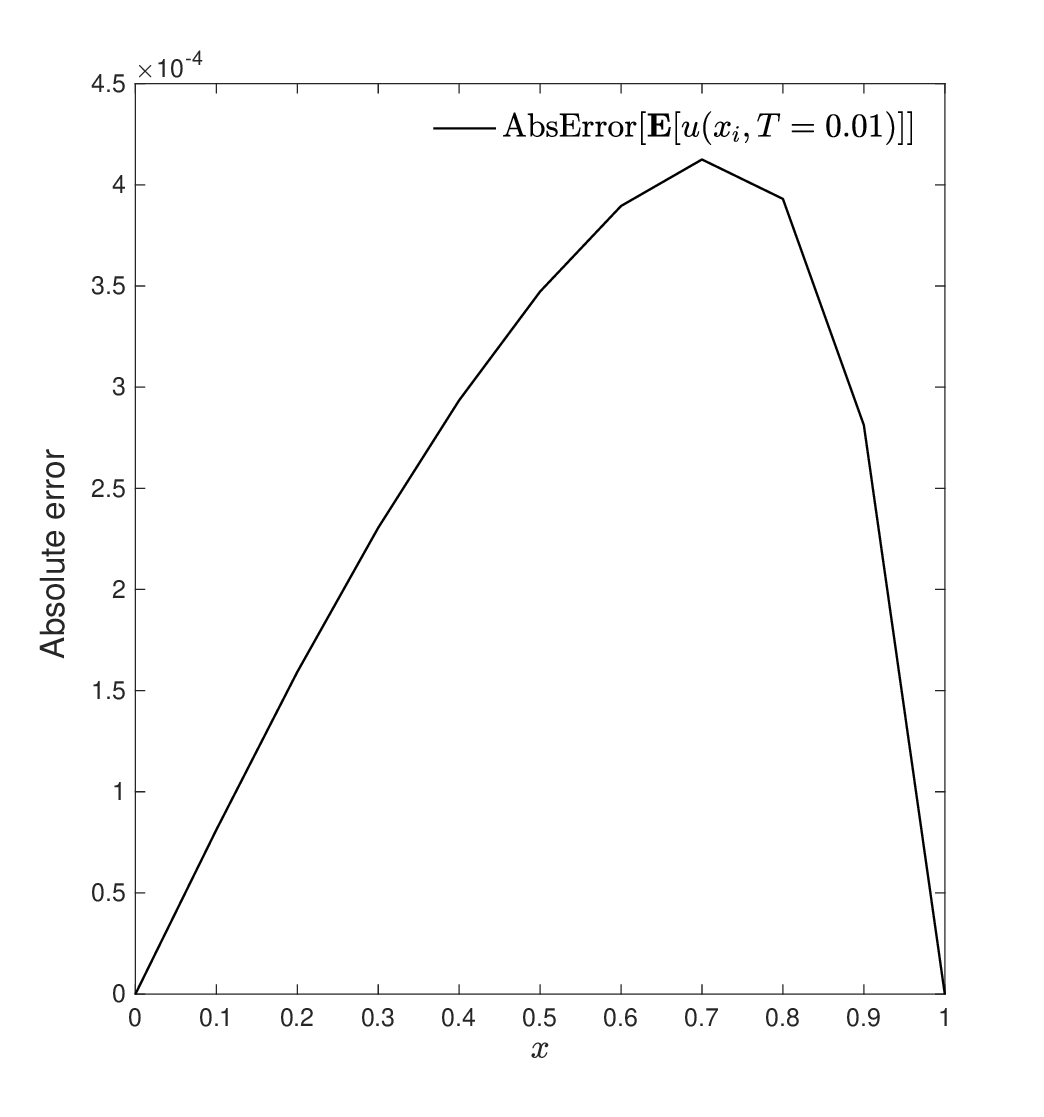}}
 \caption{\label{fig2Example} Plot (a):  Comparative graphics between the exact values of expectation of \eqref{exactSolEj}, $\mathbb{E}[u(x_i,T)]$, and the approximated expectation, $\mathbb{E}[u_i^{N_T}]$, using the random numerical scheme \eqref{EsquemaRandom}. Plot (b): Absolute error of the expectations represented in plot (a). In both graphics we have considered $T=0.01$, i.e $N_T=5$ for $k=0.002$, and $h=0.1$.}
\end{figure}

\begin{figure}[!ht]
\centering \subfigure[]%
    {\includegraphics[width=.5\textwidth]{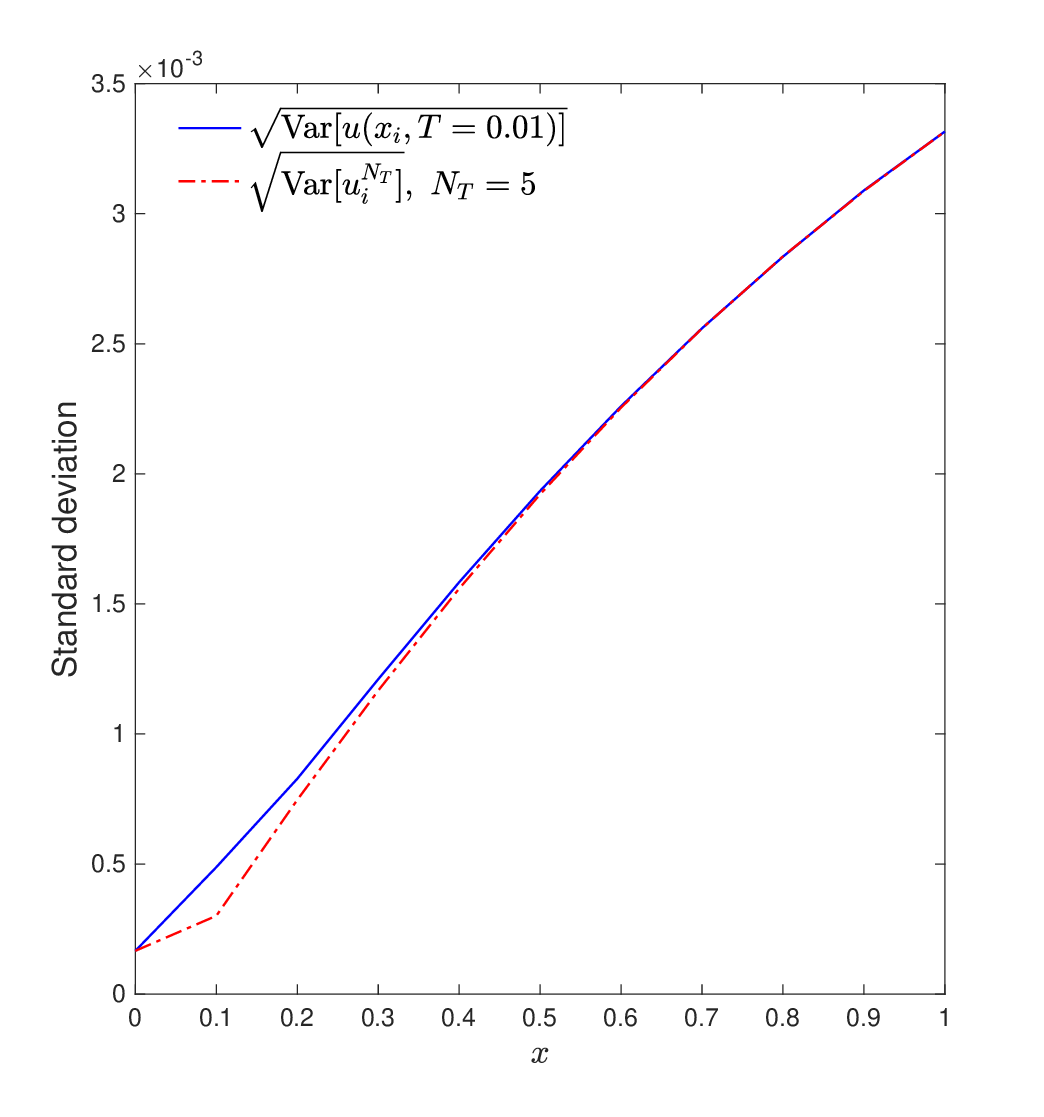}}
\!\!\!\! 
\subfigure[]
{\includegraphics[width=.5\textwidth]
{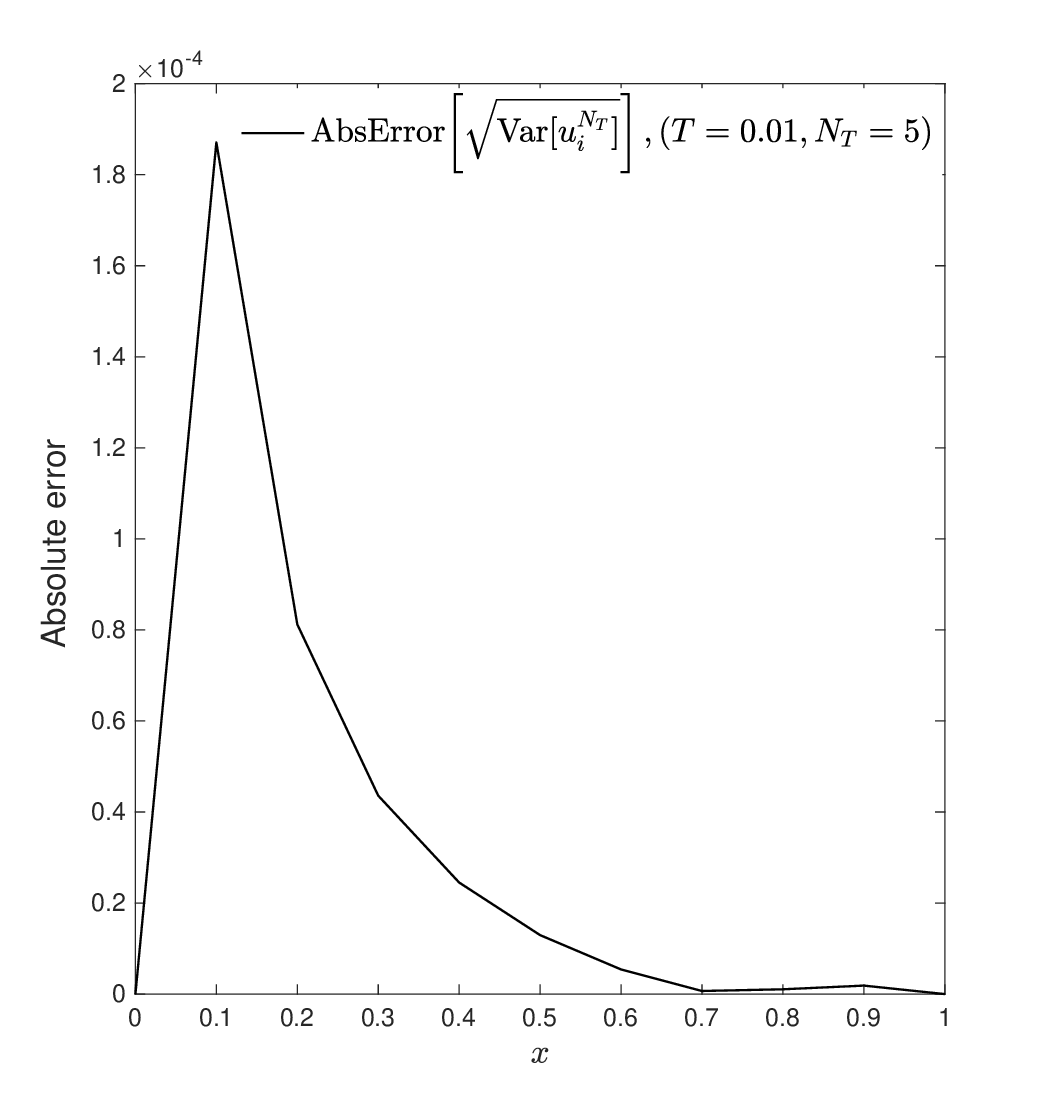}}
 \caption{\label{fig3Example} Plot (a):  Comparative graphics between the exact values of standard deviation of \eqref{exactSolEj}, $\sqrt{\mathrm{Var}[u(x_i,T)]}$, and the approximated standard deviation, $\sqrt{\mathrm{Var}\left[u_i^{N_T}\right]}$, using the random numerical scheme \eqref{EsquemaRandom}. Plot (b): Absolute error of the standard deviations represented in plot (a). In both graphics we have considered $T=0.01$, i.e $N_T=5$ for $k=0.002$, and $h=0.1$.}
\end{figure}

In Figure \ref{fig1Example} we have plotted the evolution of the exact mean (plot (a)) and the exact standard deviation (plot (b)) over the domain $(x,t)\in{ [0, 1] \times [0, 0.01]}$.
In Figures 	\ref{fig2Example} and \ref{fig3Example} we have computed the expectation and standard deviation of \eqref{exactSolEj} and check that the numerical values, obtained using Algorithm \ref{algorithm}, are close to the corresponding exact ones computing the absolute errors. Computations have been carried out using the 
software \textit{Mathematica}.

\section{Conclusions} \label{sec6}
Incorporating randomness into the mathematical models improves the quality of the approximation to real problems in the measure of the uncertainties are taken into account. This challenge involves the proof of the new intermediate results from both the analytic and numerical points of view. In this paper we use a mean square approach to random Fisher-KPP models as well as a semidiscretization technique allowing the use of previous results of systems of random ordinary differential equations. Numerical solutions are generated by constructing a random difference scheme in two steps. Firstly, we use a semidiscretization technique and then the full discretization is achieved by adapting to the random framework the ideas of the ETD method \cite{CoxMathews}. Once the random difference scheme is built, a sample approach combined with the results of \cite{DAR} recently obtained for the deterministic case permit the treatement of the stability. In spite of the high level of complexity numerical experiments illustrate the utility of the approach.

\section{Acknowledgements}

This work has been partially supported by the Ministerio de Econom\'{i}a y Competitividad grant MTM2017-89664-P.

\section*{Conflict of Interest Statement}
The authors declare that there is no conflict of interests regarding the publication of this article.

%

\end{document}